\documentclass[a4paper,10pt]{article}

\usepackage{amsmath}
\usepackage[authoryear]{natbib}
\usepackage{amsthm}
\usepackage{microtype}
\usepackage[hmargin=3cm,vmargin={3.5cm,4cm}]{geometry}
\usepackage[bitstream-charter]{mathdesign}
\usepackage[final]{showkeys}

\numberwithin{equation}{section}
\allowdisplaybreaks

\newtheorem{Theorem}{Theorem}[section]
\newtheorem{Remark}{Remark}[section]

\author{$\text{Enzo Orsingher}_1$, $\text{Federico Polito}_2$\\
	\footnotesize (1) -- Dipartimento di Scienze Statistiche,
	``Sapienza'' Universit\`a di Roma\\
	\footnotesize Piazzale Aldo Moro 5, 00185 Rome, Italy.\\
	\footnotesize Tel: +39-06-49910585, fax: +39-06-4959241\\
	\footnotesize Email address: enzo.orsingher@uniroma1.it (corresponding author)\\
	\footnotesize (2) -- Dipartimento di Matematica, Universit\`a degli studi di Roma
	``Tor Vergata''\\
	\footnotesize Via della Ricerca Scientifica 1, 00133, Rome, Italy.\\
	\footnotesize Tel: +39-06-2020568, fax: +39-06-20434631\\
	\footnotesize Email address: polito@nestor.uniroma2.it
	}

\title{The Space-Fractional Poisson Process}

\begin{document}

	\maketitle

	\begin{abstract}
		
		\noindent In this paper we introduce the space-fractional Poisson process whose state probabilities
		$p_k^\alpha(t)$, $t>0$, $\alpha \in (0,1]$, are governed by the equations
		$(\mathrm d/\mathrm dt)p_k(t) = -\lambda^\alpha (1-B)p_k^\alpha(t)$, where
		$(1-B)^\alpha$ is the fractional difference operator found in the study of time series analysis.
		We explicitly obtain the distributions $p_k^\alpha(t)$, the probability generating
		functions $G_\alpha(u,t)$, which are also expressed as distributions of the minimum
		of i.i.d.\ uniform random variables.
		The comparison with the time-fractional Poisson process is investigated and finally,
		we arrive at the more general space-time fractional Poisson process of which we give
		the explicit distribution.
		
		\medskip
		
		
		\medskip

	\end{abstract}

	\section{Introduction}
	
		Fractional Poisson processes studied so far have been obtained either by considering renewal processes
		with intertimes between events represented by Mittag--Leffler distributions \citep{mainardi,beg}
		or by replacing the time derivative in the equations governing the state probabilities with
		the fractional derivative in the sense of Caputo.
		
		In this paper we introduce a space-fractional Poisson process by means of the fractional
		difference operator
		\begin{align}
			\label{difference}
			\Delta^\alpha = (1-B)^\alpha, \qquad \alpha \in (0,1],
		\end{align}
		which often appears in the study of long memory time series \citep{tsay}.
		
		The operator \eqref{difference} implies a dependence of the state probabilities $p_k^\alpha(t)$ from
		all probabilities $p_j^\alpha(t)$, $j<k$. For $\alpha=1$ we recover the classical homogeneous
		Poisson process and the state probabilities $p_k(t)$ depend only on $p_{k-1}(t)$.
		
		For the space-fractional Poisson process we obtain the following distribution:
		\begin{align}
			\label{svetaintro}
			p_k^\alpha(t) = \text{Pr} \{ N^\alpha(t) = k \}
			= \frac{(-1)^k}{k!} \sum_{r=0}^\infty \frac{(-\lambda^\alpha t)^r}{r!} \frac{\Gamma(\alpha r+1)}{
			\Gamma(\alpha r+1-k)}, \qquad k \geq 0.
		\end{align}
		
		The distribution of the space-fractional Poisson process can be compared with that of the
		time-fractional Poisson process $N_\nu(t)$, $t>0$, $\nu \in (0,1]$:
		\begin{align}
			\label{sv}
			\text{Pr} \{ N_\nu(t) = k \} = \frac{(\lambda t^\nu)^k}{k!} \sum_{r=0}^\infty
			\frac{(r+k)!}{r!} \frac{(-\lambda t^\nu)^r}{\Gamma(\nu(k+r)+1)}, \qquad k \geq 0.
		\end{align}
		For $\alpha=\nu=1$, from \eqref{svetaintro} and \eqref{sv}, we immediately arrive at the classical
		distribution of the homogeneous Poisson process.
		
		The space-fractional Poisson process has probability generating function
		\begin{align}
			G_\alpha(u,t) = \mathbb{E}u^{N^\alpha(t)} = e^{-\lambda^\alpha(1-u)^\alpha t}, \qquad |u|\leq 1,
		\end{align}
		and can be compared with its time-fractional counterpart
		\begin{align}
			{}_\nu G(u,t) = \mathbb{E}u^{N_\nu(t)} = E_\nu\left(-\lambda(1-u)t^\nu\right), \qquad |u|\leq 1,
		\end{align}
		where
		\begin{align}
			E_\nu(x) = \sum_{r=0}^\infty \frac{x^r}{\Gamma(\nu r+1)}, \qquad \nu >0,
		\end{align}
		is the one-parameter Mittag--Leffler function.
		
		We show below that the probability generating function of the space-time fractional Poisson
		process reads
		\begin{align}
			\label{rew}
			{}_\nu G_\alpha(u,t) = E_\nu\left( -\lambda^\alpha (1-u)^\alpha t^\nu \right), \qquad |u|\leq 1,
		\end{align}
		and its distribution has the form
		\begin{align}
			p_k^{\alpha,\nu} (t) = \frac{(-1)^k}{k!} \sum_{r=0}^\infty \frac{(-\lambda^\alpha t^\nu)^r}{
			\Gamma(\nu r+1)} \frac{\Gamma(\alpha r+1)}{\Gamma(\alpha r + 1- k)}, \qquad k \geq 0, \: \alpha \in
			(0,1], \: \nu \in (0,1].
		\end{align}
		
		We also show that the space-fractional Poisson process $N^\alpha(t)$ can be regarded as a homogeneous
		Poisson process $N(t)$, subordinated to a positively skewed stable process $S^\alpha(t)$ with Laplace
		transform
		\begin{align}
			\mathbb{E} e^{-z S^\alpha(t)} = e^{-tz^\alpha}, \qquad z>0, \: t> 0.
		\end{align}
		In other words, we have the following equality in distribution
		\begin{align}
			\label{wer}
			N^\alpha(t) \overset{\text{d}}{=} N\left(S^\alpha(t)\right).
		\end{align}
		The representation \eqref{wer} is similar to the following representation of the time-fractional
		Poisson process:
		\begin{align}
			N_\nu(t) \overset{\text{d}}{=} N\left( T_{2\nu}(t) \right),
		\end{align}
		where $T_{2\nu}(t)$, $t>0$, is a process whose one-dimensional distribution is obtained by
		folding the solution to the time-fractional diffusion equation \citep{beg}.

		Finally we can note that the probability generating function \eqref{rew}, for all $u \in (0,1)$,
		can be represented as
		\begin{align}
			{}_\nu G_\alpha(u,t) = \text{Pr} \left\{ \min_{0\leq k \leq N_\nu(t)} X_k^{1/\alpha}
			\geq 1-u \right\},
		\end{align}
		where the $X_k$s are i.i.d.\ uniformly distributed random variables.
	
	\section{Construction of the space-fractional Poisson process}
	
		In this section we describe the construction of an alternative fractional generalisation of the
		classical homogeneous Poisson process. First, let us recall some basic properties.
		Let us call
		\begin{align}
			P_k(t) = \text{Pr} \{ N(t)=k \} =
			e^{-\lambda t} \frac{(\lambda t)^k}{k!}, \qquad t> 0, \: \lambda > 0, \: k \geq 0,
		\end{align}
		the state probabilities of the classical homogeneous Poisson process $N(t)$, $t>0$, of parameter
		$\lambda>0$. It is well-known that the probabilities $p_k(t)$, $k\geq 0$, solve the Cauchy problems
		\begin{align}
			\label{lena}
			\begin{cases}
				\frac{\mathrm d}{\mathrm dt} p_k(t) = -\lambda p_k(t) + \lambda p_{k-1}(t), \\
				p_k(0) =
				\begin{cases}
					0, & k>0, \\
					1, & k=0.
				\end{cases}
			\end{cases}
		\end{align}
		Starting from \eqref{lena}, some time-fractional generalisations of the homogeneous Poisson process have been
		introduced in the literature (see e.g.\ \citet{laskin, mainardi, beg}). These works are based on the substitution
		of the integer-order derivative operator appearing in \eqref{lena}
		with a fractional-order derivative operator, such as the
		Riemann--Liouville fractional derivative (as in \citet{laskin})
		or the Caputo fractional derivative (as in \citet{beg}).
		In this paper instead, we generalise the integer-order space-difference operator as follows. First, we rewrite
		equation \eqref{lena} as
		\begin{align}
			\label{lena2}
			\begin{cases}
				\frac{\mathrm d}{\mathrm dt} p_k(t) = -\lambda(1-B) p_k(t), \\
				p_k(0) =
				\begin{cases}
					0, & k>0, \\
					1, & k=0, 
				\end{cases}
			\end{cases}
		\end{align}
		where $B$ is the so-called \emph{backward shift operator} and is such that $B(p_k(t))=p_{k-1}(t)$ and
		$B^r(p_k(t)) = B^{r-1}(B(p_k(t))) = p_{k-r}(t)$. The fractional difference operator
		$\Delta^\alpha = (1-B)^\alpha$ has been widely used in time series analysis
		for constructing processes displaying long memory, such as the autoregressive fractionally integrated
		moving average process (ARFIMA). For more information on long memory processes and
		fractional differentiation, the reader can consult \citet{tsay}, page 89.
		
		Formula \eqref{lena2} can now be easily
		generalised by writing
		\begin{align}
			\label{lena3}
			\begin{cases}
				\frac{\mathrm d}{\mathrm dt} p_k^\alpha(t) = -\lambda^\alpha(1-B)^\alpha p_k^\alpha(t),
				& \alpha \in (0,1], \\
				p_k^\alpha(0) =
				\begin{cases}
					0, & k>0, \\
					1, & k=0, 
				\end{cases}
			\end{cases}
		\end{align}
		where $p_k^\alpha(t)$, $k\geq 0$, $t>0$, represents the state probabilities of a space-fractional
		homogeneous Poisson process $N^\alpha(t)$, $t>0$, i.e.\
		\begin{align}
			p_k^\alpha(t) = \text{Pr} \{ N^\alpha(t) = k \}, \qquad k \geq 0.
		\end{align}
		In turn, we have that \eqref{lena3} can also be written as
		\begin{align}
			\label{ira}
			\begin{cases}
				\frac{\mathrm d}{\mathrm dt} p_k^\alpha(t) = -\lambda^\alpha \sum_{r=0}^k \frac{\Gamma(\alpha+1)}{
				r! \Gamma(\alpha+1-r)}(-1)^r p_{k-r}^\alpha(t), & \alpha \in (0,1], \\
				p_k^\alpha(0) =
				\begin{cases}
					0, & k>0, \\
					1, & k=0.
				\end{cases}
			\end{cases}
		\end{align}
		Note that, in \eqref{ira} we considered that $p_j^\alpha(t) = 0$, $j \in \mathbb{Z}^-$.
		Equation \eqref{ira} can also be written as
		\begin{align}
			\label{bicchiere}
			\begin{cases}
				\frac{\mathrm d}{\mathrm dt} p_k^\alpha(t) = -\lambda^\alpha p_k^\alpha(t) +
				\alpha \lambda^\alpha p_{k-1}^\alpha(t) -\frac{\alpha(\alpha-1)}{2!} p_{k-2}^\alpha(t) +
				\dots + (-1)^{k+1} \frac{\alpha(\alpha-1)\dots (\alpha-k+1)}{k!} p_0^\alpha(t), \\
				p_k^\alpha(0) =
				\begin{cases}
					0, & k>0, \\
					1, & k=0.
				\end{cases}
			\end{cases}
		\end{align}
		By applying the reflection property of the gamma function $\Gamma(z)\Gamma(1-z) = \pi /\sin (\pi z)$
		for $z=r-\alpha$, we have also that
		\begin{align}
			\label{iranew}
			\begin{cases}
				\frac{\mathrm d}{\mathrm dt} p_k^\alpha(t) = -\lambda^\alpha p_k^\alpha(t) + \alpha \lambda^\alpha
				p_{k-1}^\alpha(t) + \frac{\lambda^\alpha \sin(\pi \alpha)}{\pi} \sum_{r=2}^k
				B(\alpha +1,r-\alpha) p_{k-r}^\alpha(t), \\
				p_k^\alpha(0) =
				\begin{cases}
					0, & k>0, \\
					1, & k=0,
				\end{cases}
			\end{cases}
		\end{align}
		where the sum is considered equal to zero when $k<2$ and $B(x,y)$ is the beta function.
		From \eqref{bicchiere} and \eqref{iranew} we see that for $\alpha=1$ we retrieve equation \eqref{ira}
		of the homogeneous Poisson process.
		
		An example of process whose state probabilities $\tilde{p}_k(t)$, depend on all $\tilde{p}_j(t)$, $j < k$,
		is the iterated Poisson process $\tilde{N}(t) = N_1(N_2(t))$, where $N_1(t)$ and $N_2(t)$ are
		independent homogeneous Poisson processes. The process $\tilde{N}(t)$ is analysed in \citet{comp}.

		\begin{Theorem}
			Let $N^\alpha(t)$, $t>0$, be a space-fractional homogeneous Poisson process of parameter $\lambda>0$ and
			let $G_\alpha(u,t) = \mathbb{E}u^{N^\alpha(t)}$, $|u|\leq 1$, $\alpha \in (0,1]$,
			be its probability generating function.
			The Cauchy problem satisfied by $G_\alpha(u,t)$ is
			\begin{align}
				\label{tomsk}
				\begin{cases}
					\frac{\partial}{\partial t} G_\alpha(u,t) = -\lambda^\alpha
					G_\alpha(u,t) (1-u)^\alpha, & |u| \leq 1, \\
					G_\alpha(u,0) = 1,
				\end{cases}
			\end{align}
			with solution
			\begin{align}
				\label{tomsk2}
				G_\alpha(u,t) = e^{-\lambda^\alpha t(1-u)^\alpha}, \qquad |u| \leq 1,
			\end{align}
			that is, the probability generating function of a discrete stable distribution
			(see \citet{devroye}, page 349).
			
			\begin{proof}
				Starting from \eqref{ira}, we have that
				\begin{align}
					\frac{\partial}{\partial t} G_\alpha(u,t) & = -\lambda^\alpha \Gamma(\alpha +1)
					\sum_{r=0}^\infty \sum_{k=r}^\infty u^k \frac{(-1)^r}{r! \Gamma(\alpha+1-r)} p_{k-r}^\alpha (t) \\
					& = -\lambda^\alpha \Gamma(\alpha+1) \sum_{r=0}^\infty \sum_{k=0}^\infty \frac{u^{k+r}(-1)^r}{
					r! \Gamma(\alpha+1-r)}p_k^\alpha(t) \notag \\
					& = -\lambda^\alpha
					\Gamma(\alpha+1) G_\alpha(u,t) \sum_{r=0}^\infty \frac{(-1)^r}{r!\Gamma(\alpha+1-r)}
					\notag \\
					& = -\lambda^\alpha G_\alpha(u,t) (1-u)^\alpha, \notag
				\end{align}
				thus obtaining formula \eqref{tomsk}. It immediately follows that
				\begin{align}
					G_\alpha(u,t) = e^{-\lambda^\alpha t (1-u)^\alpha}, \qquad |u|\leq 1.
				\end{align}
			\end{proof}
		\end{Theorem}
		
		\begin{Remark}
			Note that, for $\alpha=1$,
			formula \eqref{tomsk2} reduces to the probability generating function of the classical
			homogeneous Poisson process. Furthermore, from \eqref{tomsk2} we have that
			$\mathbb{E}[N^\alpha(t)]^j = \infty$, $j \in \mathbb{N}$, $\alpha \in (0,1)$.
		\end{Remark}
		
		\begin{Remark}
			Let $X_k$, $k=1,\dots$, be i.i.d.\ Uniform$[0,1]$ random variables, then
			\begin{align}
				G_\alpha(u,t) = e^{-\lambda^\alpha
				t(1-u)^\alpha} = \text{Pr} \left\{ \min_{0\leq k \leq N(t)} X_k^{1/\alpha}
				\geq 1-u \right\}, \qquad u \in (0,1),
			\end{align}
			where $N(t)$, $t>0$, is a classical homogeneous Poisson process of rate $\lambda^\alpha$
			with the assumption that $\min ( X_k^{1/\alpha} ) = 1$ when $N(t)=0$.
		\end{Remark}
		
		\begin{Theorem}
			The discrete stable state probabilities of a space-fractional homogeneous Poisson process
			$N^\alpha(t)$, $t>0$, can be written as
			\begin{align}
				\label{sveta}
				p_k^\alpha(t) & = \text{Pr} \{ N^\alpha(t) = k \} \\
				& = \frac{(-1)^k}{k!} \sum_{r=0}^\infty \frac{(-\lambda^\alpha t)^r}{r!} \frac{\Gamma(\alpha r+1)}{
				\Gamma(\alpha r+1-k)} \notag \\
				& = \frac{(-1)^k}{k!} \: {}_1\psi_1 \left[ \left.
				\begin{array}{l}
					(1,\alpha) \\
					(1-k,\alpha)
				\end{array}
				\right| -\lambda^\alpha t \right], \qquad t>0, \: k \geq 0, \notag
			\end{align}
			where ${}_h\psi_j(z)$ is the generalised Wright function (see \citet{kilbas}, page 56,
			formula (1.11.14)).
			
			\begin{proof}
				By expanding the probability generating function \eqref{tomsk2} we have that
				\begin{align}
					G_\alpha(u,t) & = e^{-\lambda^\alpha t (1-u)^\alpha} \\
					& = \sum_{r=0}^\infty \frac{[-\lambda^\alpha t(1-u)^\alpha]^r}{r!} \notag \\
					& = \sum_{r=0}^\infty \frac{(-\lambda^\alpha t)^r}{r!} \sum_{m=0}^\infty
					\frac{(-u)^m \Gamma(\alpha r+1)}{m! \Gamma(\alpha r +1 -m)} \notag \\
					& = \sum_{m=0}^\infty u^m \frac{(-1)^m}{m!} \sum_{r=0}^\infty
					\frac{(-\lambda^\alpha t)^r}{r!} \frac{\Gamma(\alpha r+1)}{\Gamma(\alpha r+1-m)}. \notag
				\end{align}
				From this, formula \eqref{sveta} immediately follows.
			\end{proof}
		\end{Theorem}
		
		\begin{Remark}
			Note that the discrete stable distribution \eqref{sveta} (which for $\alpha=1$ reduces to the Poisson
			distribution) can also be written as
			\begin{align}
				p_k^\alpha(t) & = \frac{(-1)^k}{k!} \sum_{r=0}^\infty \frac{(-\lambda^\alpha t)^r}{r!}
				\frac{\Gamma(\alpha r+1)}{\Gamma(\alpha r+1-k)} \\
				& = \frac{(-1)^k}{k!} \int_0^\infty e^{-w} \sum_{r=0}^\infty \frac{(-\lambda^\alpha t w^\alpha)^r}{
				r! \Gamma(\alpha r+1-k)}. \notag
			\end{align}
		\end{Remark}

		\begin{Theorem}
			Let $S^\gamma(t)$, $t>0$, $\gamma \in (0,1)$, be a $\gamma$-stable subordinator, that is a positively skewed
			stable process such that
			\begin{align}
				\mathbb{E}e^{-zS^\gamma(t)} = e^{-tz^\gamma}, \qquad z>0, \: t>0,
			\end{align}
			and with transition function $q_\gamma(s,t)$.
			For a space-fractional Poisson process $N^\alpha(t)$, $t>0$, $\alpha \in (0,1]$, with
			rate $\lambda>0$, the following
			representation holds in distribution:
			\begin{align}
				\label{volpe}
				N^\alpha(S^\gamma(t)) \overset{\text{d}}{=} N^{\alpha \gamma}(t).
			\end{align}
			
			\begin{proof}
				In order to prove the representation \eqref{volpe} it is sufficient to observe that
				\begin{align}
					\int_0^\infty G_\alpha(u,s) q_\gamma(s,t) \mathrm ds
					= \int_0^\infty e^{-\lambda^\alpha s(1-u)^\alpha} q_\gamma(s,t) \mathrm ds
					= e^{-\lambda^{\alpha \gamma} t (1-u)^{\alpha \gamma}}.
				\end{align}
			\end{proof}
		\end{Theorem}
		
		\begin{Remark}
			Note that, when $\alpha=1$, formula \eqref{volpe} reduces to
			\begin{align}
				N(S^\gamma(t)) \overset{\text{d}}{=} N^\gamma(t),
			\end{align}
			and this reveals a second possible way of constructing a space-fractional Poisson process.
		\end{Remark}
		
		Consider now the first-passage time at $k$ of the space-fractional Poisson process
		\begin{align}
			\tau_k^\alpha (t) = \inf \{ t \colon N^\alpha(t) = k \}, \qquad k \geq 0. 
		\end{align}
		Since $\text{Pr} \{ \tau_k^\alpha < t \} = \text{Pr} \{ N^\alpha(t) \geq k \}$, we have that
		\begin{align}
			\text{Pr} \{ \tau_k^\alpha < t \} = \sum_{m=k}^\infty \frac{(-1)^m}{m!}
			\sum_{r=0}^\infty \frac{(-\lambda^\alpha t)^r}{r!} \frac{\Gamma(\alpha r+1)}{\Gamma(\alpha r +1-m)}.
		\end{align}
		Hence
		\begin{align}
			\text{Pr} \{ \tau_k^\alpha \in \mathrm ds \} = \sum_{m=k}^\infty \frac{(-1)^m}{m!}
			\sum_{r=1}^\infty \frac{(-\lambda^\alpha)^r t^{r-1}}{(r-1)!} \frac{\Gamma(\alpha r+1)}{\Gamma(\alpha r +1-m)}.
		\end{align}
		Note that for $\alpha=1$ we obtain the classical Erlang process. First we have
		\begin{align}
			\text{Pr} \{ \tau_k^1 < t \} & = \sum_{m=k}^\infty \frac{(-1)^m}{m!}
			\sum_{r=m}^\infty \frac{(-\lambda t)^r}{(r-m)!} \\
			& = e^{-\lambda t} \sum_{m=k}^\infty \frac{(\lambda t)^m}{m!}, \notag
		\end{align}
		and therefore
		\begin{align}
			\text{Pr} \{ \tau_k^1 \in \mathrm ds \} & = -\lambda e^{-\lambda t} \sum_{m=k}^\infty
			\frac{(\lambda t)^m}{m!} + \lambda e^{-\lambda t} \sum_{m=k}^\infty \frac{(\lambda t)^{m-1}}{(m-1)!} \\
			& = -\lambda e^{-\lambda t} \sum_{m=k}^\infty
			\frac{(\lambda t)^m}{m!} + \lambda e^{-\lambda t} \sum_{m=k-1}^\infty \frac{(\lambda t)^m}{m!} \notag \\
			& = \lambda e^{-\lambda t} \frac{(\lambda t)^{k-1}}{(k-1)!}. \notag
		\end{align}
		
		\begin{Remark}
			Note that, with a little effort, fractionality can be introduced also in time.
			In this case, for example, the fractional differential equation governing the
			probability generating function is
			\begin{align}
				\label{fractomsk}
				\begin{cases}
					\frac{\partial^\nu}{\partial t^\nu} {}_\nu G_\alpha(u,t) = -\lambda^\alpha
					{}_\nu G_\alpha(u,t) (1-u)^\alpha, & |u| \leq 1, \: \nu \in (0,1], \: \alpha \in (0,1], \\
					{}_\nu G_\alpha(u,0) = 1,
				\end{cases}
			\end{align}
			where $\partial^\nu/\partial t^\nu$ is the Caputo fractional derivative operator (see \citet{kilbas}).
			By means of Laplace transforms, some simple manipulations lead to
			\begin{align}
				{}_\nu G_\alpha(u,t) = E_\nu(-\lambda^\alpha t^\nu (1-u)^\alpha), \qquad |u| \leq 1,
			\end{align}
			where $E_\nu(x)$ is the Mittag--Leffler function \citep{kilbas}. In turn, by expanding the
			above probability generating function we have that
			\begin{align}
				p_k^{\alpha,\nu} (t) = \frac{(-1)^k}{k!} \sum_{r=0}^\infty \frac{(-\lambda^\alpha t^\nu)^r}{
				\Gamma(\nu r+1)} \frac{\Gamma(\alpha r+1)}{\Gamma(\alpha r + 1- k)}, \qquad k \geq 0, \: \alpha \in
				(0,1], \: \nu \in (0,1].
			\end{align}
			When $\alpha=1$ these probabilities easily reduce to those of a fractional Poisson process
			(see \citet{beg}):
			\begin{align}
				p_k^{1,\nu}(t) & = \frac{(-1)^k}{k!} \sum_{r=0}^\infty \frac{(-\lambda t^\nu)^r}{\Gamma(\nu r+1)}
				\frac{\Gamma(r+1)}{\Gamma(r-k+1)} \\
				& = \frac{(-1)^k}{k!} \sum_{r=k}^\infty \frac{(-\lambda t^\nu)^r}{\Gamma(\nu r+1)}
				\frac{r!}{(r-k)!} \notag \\
				& = \sum_{r=k}^\infty (-1)^{r-k} \binom{r}{k} \frac{(\lambda t^\nu)^r}{\Gamma(\nu r+1)}. \notag 
			\end{align}
			Moreover, let $X_k$, $k=1,\dots$, be i.i.d.\ Uniform$[0,1]$ random variables and $N_\nu(t)$, $t>0$, be a
			homogeneous time-fractional Poisson process of rate $\lambda^\alpha$ with the assumption that
			$\min ( X_k^{1/\alpha} ) = 1$ when $N_\nu(t) = 0$. The probability generating function
			${}_\nu G_\alpha(u,t)$, for $u \in (0,1)$, can be written as
			\begin{align}
				\label{prpr}
				{}_\nu G_\alpha(u,t) & = E_\nu(-\lambda^\alpha t^\nu (1-u)^\alpha) \\
				& = \sum_{r=0}^\infty (-1)^r \frac{(\lambda^\alpha t^\nu)^r(1-u)^{\alpha r}}{\Gamma(\nu r+1)} \notag \\
				& = \sum_{r=0}^{\infty} (-1)^r \frac{(\lambda^\alpha t^\nu)^r}{\Gamma(\nu r+1)}
				\sum_{k=0}^r (-1)^k \binom{r}{k} \left[1-(1-u)^\alpha\right]^k \notag \\
				& = \sum_{k=0}^\infty \left[ 1-(1-u)^\alpha \right]^k \sum_{r=k}^\infty (-1)^{r-k}
				\binom{r}{k} \frac{(\lambda^\alpha t^\nu)^r}{\Gamma(\nu r+1)} \notag \\
				& = \sum_{k=0}^\infty \left[ \text{Pr} \left( X_k^{1/\alpha} \geq 1-u \right) \right]^k
				\text{Pr} \left\{ N_\nu(t) = k \right\} \notag \\
				& = \text{Pr} \left\{ \min_{0\leq k \leq N_\nu(t)} X_k^{1/\alpha}
				\geq 1-u \right\}. \notag
			\end{align}
			Formula \eqref{prpr} shows that the contribution of the space-fractionality affects only the
			uniform random variables $X_k^{1/\alpha}$ while the time-fractionality only the driving process
			$N_\nu(t)$.
		\end{Remark}

	\bibliographystyle{abbrvnat}
	\bibliography{space}
	\nocite{*}

\end{document}